\documentclass[12pt]{article}
\usepackage[T1]{fontenc}
\usepackage[ansinew]{inputenc}
\usepackage{amsmath}
\usepackage{graphicx}
\title{Manifolds with weighted Poincar\'e inequality and uniqueness of minimal hypersurfaces\thanks{Both authors are  partially supported by CNPq and Faperj  of Brazil.}}
\author{Xu Cheng and Detang Zhou}
\usepackage{amsmath, amssymb}
\usepackage{amsthm}
\makeatletter

\@addtoreset{equation}{section} \makeatother
\newtheorem{thm}{Theorem}[section]
 
\newtheorem{cor}{Corollary}[section]

\begin{document}
\date{}
\maketitle
 \begin{abstract}
In this paper, we obtain  results on rigidity of complete Riemannian
manifolds with weighted Poincar\'e inequality. As an application, we
prove that if $M$ is a complete $\frac{n-2}{n}$-stable minimal
hypersurface in $\mathbb{R}^{n+1}$ with $n\geq 3$ and has bounded
norm of the second fundamental form, then $M$ must either have only
one end or be a catenoid.

\end{abstract}

\baselineskip=14pt
\section {Introduction}

In this paper, we will discuss complete Riemannian manifolds with
weighted Poincar\'e inequality and  minimal hypersurfaces with
$\delta$-stability in the Euclidian space $\mathbb{R}^{n+1}$ with
$n\geq 3$. We first recall some backgrounds.

Let $M$ be an $n$-dimensional Riemannian manifold. Given a
Schr\"odinger operator $L=\Delta + q(x), x\in M$ on $M$, we consider
the  eigenvalue problem on a compact subdomain $D\subset M$:
$$
\left\{
  \begin{array}{ll}
    Lf+\lambda f=0, & \hbox{\text{ in }  $D$;} \\
    f|_{\partial D}=0. & \hbox{}
  \end{array}
\right.
$$
It has discrete spectrum and the number of negative eigenvalues is
finite.  The (Morse) index of L on $M$ is defined as the supremum,
over compact domains of $M$, of the number of negative eigenvalues
(counted with multiplicity) of $L$ with Dirichlet boundary
condition.

If $M$ is a complete connected immersed  minimal hypersurface in
$\mathbb{R}^{n+1}, n\geq 2$ and if $L$ is the Jacobi operator
$L=\Delta +|A|^2$, where $|A|$ denotes the norm of the second
fundamental form of $M$, then the index of $L$ is said to be the
(Morse) index of $M$. $M$ is said to be stable if the index of $M$
is $0$, which is equivalent to say that, for all  compactly
supported piecewisely smooth function $f\in C_o^\infty(M)$,
\begin{equation}\label{stability1}
   \int_M\left(|\nabla f|^2-|A|^2 f^2 \right)\ge 0.
\end{equation}

It is known that a complete stable minimal surface in $\mathbb{R}^3$
is plane, which was proved by do Carmo and Peng, and Fischer-Cobrie
and Schoen independently (\cite{dCP}, \cite{FS}); and
 that only index one complete minimal surfaces in $\mathbb{R}^3$ are
the catenoid and Enneper surface, which was proved by Lopez and Ros
\cite{LR}.

 While it is unknown that a complete stable minimal hypersurface
  in $\mathbb{R}^{n+1}$ is a hyperplane when $n\leq 7$, Cao, Shen and
 Zhu \cite{CSZ} proved that a complete stable minimal hypersurface in
 $\mathbb{R}^{n+1} $ must have only one end for all dimension  $n\geq 3$.
 Tam and Zhou \cite{TZ} recently showed that an ($n$-diemensional) catenoid in the Euclidean space
$\mathbb{R}^{n+1}$ with  $n\geq 3$ has index $1$ .

Now let us assume $L=\Delta +\delta|A|^2$ on minimal hypersurface
$M$ in $\mathbb{R}^{n+1}$ for some number $0<\delta\leq 1$ . We may
similarly define that $M$ is {\it $\delta$-stable} if
\begin{equation}\label{stability2}
   \int_M\left(|\nabla f|^2-\delta|A|^2 f^2 \right)\ge 0
\end{equation}
 for all $f\in C_o^\infty(M)$.

Obviously, given $\delta_1>\delta_2$, $\delta_1$-stable implies
$\delta_2$-stable. So $M$ is stable implies that $M$ is
$\delta$-stable for all $0<\delta\leq 1$. Hyperplane  is
$\delta$-stable for all $0<\delta\leq 1$.

There are some work on  $\delta$-stable minimal hypersurfacs. Kawai
\cite{K} proved a $\delta$-stable, $\delta>\frac{1}{8}$ complete
minimal surface in $\mathbb{R}^3$ must be plane. Recently, Meeks,
P\'erez and Ros \cite{MPR} showed that any complete embedded
$\delta$-stable minimal surface in $\mathbb{R}^3$ with finite genus
is flat. In the case of higher dimension $n\geq 3$, we have,
directly from the argument in \cite{CSZ}, that the result of Cao,
Shen and Zhu also holds for $\frac{n-1}{n}$-stable. Recently Tam and
Zhou \cite{TZ} showed that a catenoid in $\mathbb{R}^{n+1}$ is
$\frac{n-2}{n}$-stable. Also they proved that if $M$ is an
$\frac{n-2}{n}$-stable complete immersed minimal hypersurface in
$\mathbb{R}^{n+1}$ and if
\begin{equation*}
  \lim_{R\to +\infty} \frac
1{R^2}\int_{B(2R)\setminus B(R)}|A|^{\frac{2(n-2)}{n}}=0,
\end{equation*}
then  $M$ is either a hyperplane or a catenoid.

In this paper, we prove that if  an $\frac{n-2}{n}$-stable complete
minimal hypersurface in $\mathbb{R}^{n+1}$ with $ n\geq 3$ and the
norm of its
 second fundamental form satisfies some  growth condition, then
 it either has only one end or is a catenoid. More precisely,
we show

 \begin{thm}\label{thm01}
Let $M$ be an $\frac{n-2}{n}$-stable complete minimal hypersurface
in $\mathbb{R}^{n+1}$ for $ n\geq 3$ and the norm of its
 second fundamental form satisfies
\begin{align}
    &\lim_{R\to +\infty}\sup_{B(R)}|A|
/R^{\frac{n-3}{2}}=0 \quad \text{for}\quad n\geq 4; \notag\\
&\lim_{R\to +\infty}\sup_{B(R)}|A|/{\ln R}=0 \quad \text{for}\quad
n=3,
\end{align} then $M$ either has only one end or is a
catenoid.
\end{thm}

From Theorem \ref{thm01}, we have the following result:
\begin{cor}\label{bound} Let $M$ be an $\frac{n-2}{n}$-stable complete minimal hypersurface in
$\mathbb{R}^{n+1}, n\geq 3$ with at least two ends. If it has
bounded norm of the second fundamental form, then $M$ must be a
catenoid.
\end{cor}

 Our results for minimal hypersurfaces   rely on the study on complete manifolds with
 weighted Poincar\'e inequality which is of independent  interest.

 Let $M$ be a complete Riemannian manifold. Recall from \cite{LW3}
that a complete Riemannian manifold $(M,ds^2)$ is said to satisfy a
weighted Poincar\'e inequality with nonnegative weighted function
$\rho$ if the inequality
$$\int_M|\nabla f|^2\ge \int_M\rho f^2$$
holds for all  compactly supported piecewisely smooth function $f\in
C_o^{+\infty}(M)$.

Further, $M$ is said to satisfy property $(\mathcal{P}_{\rho})$ for
non-zero nonnegative weight function $\rho(x)$ if,

(1) $M$ satisfies a weighted Poincar\'e inequality with $\rho$; and

(2) the conformal metric $\rho ds^2$ is complete.

 In \cite{LW3}, Li and Wang studied  complete manifolds
satisfying
 property $(\mathcal{P}_\rho)$ and obtained many theorems on rigidity. Later the first author \cite{C} discussed
complete manifolds with Poincar\'e inequality and obtain results on
the uniqueness of ends which can be applied to study stable minimal
hypersurfaces in a Riemannian manifold.   In this paper, we
generalize one result of Li and Wang (\cite{LW3}, Theorem 5.2) to
the following

\begin{thm}\label{thm-r01}Let $M$ be a  complete $n$-dimensional $(n\geq 3)$ Riemannian manifold with property $(\mathcal{P}_{\rho})$ for some
nonzero weight function $\rho$. Suppose the Ricci curvature of $M$
has the lower bound
\begin{equation*}
    Ric_M(x)\ge -(n-1) \tau(x), \quad x\in M,
\end{equation*}
where $\tau(x)$ satisfies  Poincar\'e inequality
$$\int_M|\nabla f|^2\ge (n-2)\int_M\tau f^2, \textrm{ for all
  }f\in C_o^{+\infty}(M).$$ If $\rho$ and $\tau $ satisfy the
growth condition
\begin{align}\label{eq}
    &\liminf_{R\to +\infty}S(R)
e^{-\frac{(n-3)}{n-2}R}=0\quad \text{for}\quad n\geq 4 \notag\\
&\liminf_{R\to +\infty}S(R){R}^{-1}=0 \quad \text{for}\quad  n=3,
\end{align}
where
\[
    S(R) =\sup_{x\in B_\rho(R)}(\sqrt{\rho(x)},\sqrt{\tau(x)}),
  \]
then either
\item[1)] $M$ has only one nonparabolic end; or
\item[2)] $M$ has two nonparabolic ends and is given by $M=\mathbb{R}\times N$
with the warped product metric
$$ds_M^2=dt^2+\eta^2(t)ds^2_N,$$
for some positive function $\eta(t)$ and some compact manifold $N$.
Moreover, $\tau(t)$ is a function of $t$ alone satisfying
$$(n-2)\eta''\eta^{-1}=\tau.$$

\end{thm}
If we choose $\tau=\frac{1}{n-2}\rho$ in Theorem \ref{thm-r01}, it
is just Theorem 5.2 of \cite{LW3}.  In the case of minimal
hypersurfaces, we couldn't find any weighed function $\rho$ in a
Poincar\'e inequality, which satisfies both the completeness of the
metric $\rho ds^2$ and the lower bound estimate of Ricci curvature
of $M$. Hence we couldn't apply the theorem of Li and Wang. Instead,
our theorem \ref{thm-r01} is suitable to our minimal case (see
Theorem \ref{thm01}).

 The work of Li and Wang on
complete manifolds satisfying  weighted Poincar\'e inequality is a
generalization of their one on complete manifolds with positive
spectrum (\cite{LW1} and \cite{LW2}. See \cite{LW3} and the
references therein). Let $\lambda_1(M)$ be the largest lower bound
of the spectrum of the Laplacian with respect to the metric of $M$.
Theorem \ref{thm-r01} implies the following result.
\begin{cor}
\label{cor-r01}Let $M$ be a  complete $n$-dimensional $(n\geq 3)$
Riemannian manifold with positive spectrum (i.e. $\lambda_1(M)>0)$.
Suppose the Ricci curvature of $M$ has the lower bound
\begin{equation*}
    Ric_M(x)\ge -(n-1) \tau(x), \quad x\in M,
\end{equation*}
where $\tau(x)$ satisfies  Poincar\'e inequality
$$\int_M|\nabla f|^2\ge (n-2)\int_M\tau f^2, \textrm{ for all
  }f\in C_o^{+\infty}(M).$$ If  $\tau $ satisfies the
growth condition
\begin{align}\label{eq}
    &\liminf_{R\to +\infty}(\sup_{x\in
B(R)}\tau(x))
e^{-\frac{2(n-3)}{n-2}R}=0\quad \text{for}\quad n\geq 4 \notag\\
&\liminf_{R\to +\infty}(\sup_{x\in B(R)}\tau(x)){R}^{-2}=0 \quad
\text{for}\quad  n=3,
\end{align}

then either
\item[1)] $M$ has only one nonparabolic end; or
\item[2)] $M$ has two nonparabolic ends and is given by $M=\mathbb{R}\times N$
with the warped product metric
$$ds_M^2=dt^2+\eta^2(t)ds^2_N,$$
for some positive function $\eta(t)$ and some compact manifold $N$.
Moreover, $\tau(t)$ is a function of $t$ alone satisfying
$$(n-2)\eta''\eta^{-1}=\tau.$$
\end{cor}

 This corollary
generalizes Theorem 2.1 in \cite{LW1} (just choose
$\tau(x)=\frac{\lambda_1(M)}{n-2}$ and use the fact a nonparabolic
end with $\lambda_1(M)>0$ has infinite volume).

Throughout this paper, all manifolds are assumed to be oriented.
\bigskip

{\bf Acknowledgement:} The authors would like to thank  Jiaping Wang
for very helpful discussions during their  visit to the Chinese Hong
Kong University.

\section{Rigidity of complete manifolds}

In this section, we will consider the  structure of a complete
manifold $M$ with property $(\mathcal{P}_\rho)$. Since we follow the
argument of Li and Wang (\cite{LW3},  Theorem 5.2) with some changes
of techniques in the proof of our Theorem \ref{thm-r01}, we
recommend \cite{LW3} as a complement when necessary.

Let $d(x,y)$ and $d_\rho(x,y)$ denote the distance between $x$ and
$y$ with respect to  $ds^2$ and  $\rho^2ds^2$ respectively.
$B(x,R)=\{y\in M; d(x,y)<R\}$ and $B_\rho(x,R)=\{y\in M;
d_\rho(x,y)<R\}$.  For a fixed point $p\in M$, we denote $r(x)$ and
$r_{\rho}(x)$ the distance function with respect to metric $ds^2$
and conformal metric $\rho ds^2$ from $p$ respectively. $B(R)=\{x\in
M; r(x)< R\}$ and $B_{\rho}(R)=\{x\in M; r_{\rho}(x)< R\}$.

We need the following construction of harmonic functions (see
\cite{LW3}, $\S 5$).

 Suppose $M$ has at least two nonparabolic ends $E_1$ and $E_2$. A theory of Li and Tam \cite{LT} asserts that
one can get a nonconstant bounded harmonic function $f$ with finite
Dirichlet integral by taking a convergent subsequence of the
harmonic functions $f_R$ as $R\to +\infty$, satisfying
 \[ \Delta f_R=0 \qquad \textrm{ on }B(R),
 \]
 with boundary conditions
 \begin{equation*}
    \left\{
       \begin{array}{ll}
         f_R=1, & \hbox{on $\partial B(R)\cap E_1$;} \\
         f_R=0, & \hbox{on $\partial B(R)\setminus E_1$.}
       \end{array}
     \right.
 \end{equation*}
 It follows from the maximum principle that $0\le f_R\le 1$ for all $R$ and hence $
 0\le f\le 1. $

Now we prove Theorem \ref{thm-r01}.

\begin{proof} If $M$ has at least two nonparabolic ends, then  there exists a
bounded harmonic function $f$ with finite Dirichlet integral
constructed as above. We may assume that  $\inf f=0$ and $\sup f=1$.

Then the  Bochner formula  and the lower bound of the Ricci
curvature imply (cf. \cite{LW3}, Lemm 4.1)
\begin{equation}\label{boch}
    \Delta|\nabla f|\ge -(n-1)\tau|\nabla
    f|+\frac1{n-1}\frac{|\nabla|\nabla f||^2}{|\nabla f|}.
\end{equation}
Set  $\alpha=\frac{n-2}{n-1}$ and $g=|\nabla f|^{\alpha}$.
(\ref{boch}) implies
\begin{equation}\label{eq7}
    \begin{split}
       \Delta g& =\alpha(\alpha-1)|\nabla f|^{\alpha-2}|\nabla|\nabla f||^2+\alpha |\nabla f|^{\alpha -1}\Delta |\nabla f| \\
         &\ge -(n-2)\tau g.
     \end{split}
\end{equation}

 We will show
 inequality (\ref{eq7})
is actually an equality. For any  nonnegative compactly supported
piecewisely smooth function $\phi$ on $M$, we have
\begin{align}\label{eqn13}
  &\int_M\phi^2g(\Delta g +(n-2)\tau g)\notag\\
&\quad = -2\int_M\phi g\langle    \nabla g,\nabla \phi\rangle  -\int_M\phi^2|\nabla g|^2 +\int_M(n-2)\tau(\phi g)^2 \notag \\
&\quad \le  -2\int_M\phi g\langle    \nabla g,\nabla \phi\rangle  -\int_M\phi^2|\nabla g|^2 +\int_M|\nabla(\phi g)|^2  \\
 &\quad = \int_M|\nabla\phi|^2|\nabla f|^{\frac{2(n-2)}{n-1}}= \int_M|\nabla \phi|^2g^2.\notag
\end{align}
The inequality in (\ref{eqn13}) holds since $\tau$ satisfies
Poincar\'e inequality.

Choose $\phi=\psi\chi$, where $\psi$ and $\chi$ denote two
piecewisely smooth compactly supported functions on $M$ to be chosen
later. Then
\begin{equation}\label{eqn14}
    \int_M|\nabla \phi|^2g^2\le 2\int_M |\nabla \psi|^2\chi^2|\nabla
    f|^\frac{2(n-2)}{n-1}+2\int_M |\nabla \chi|^2\psi^2|\nabla
    f|^\frac{2(n-2)}{n-1}.
\end{equation}

We first consider the case of $n\geq 4$. For $R>1$ we let $\psi(x)$
be a function depending on the $\rho$-distance:
\begin{equation*}
    \psi(x)=\left\{
           \begin{array}{ll}
             1 & \hbox{on $B_{\rho}(R-1)$,} \\
             R-r_{\rho} & \hbox{on $B_{\rho}(R)\setminus B_{\rho}(R-1)$,} \\
             0 & \hbox{on $M\setminus B_{\rho}(R)$.}
           \end{array}
         \right.
\end{equation*}

 For $\sigma\in(0,1)$ and
$\epsilon\in (0,\frac12)$, we define $\chi$  on the level sets of
$f$:
\begin{equation*}
    \chi(x)=\left\{
             \begin{array}{ll}
               0 & \hbox{ on $\mathcal{L}(0,\sigma\epsilon)\cup \mathcal{L}(1-\sigma\epsilon,1)$,} \\
               (\epsilon-\sigma\epsilon)^{-1}(f-\sigma\epsilon)  & \hbox{on $\mathcal{L}(\sigma\epsilon,\epsilon)\cap(M\setminus E_1)$,} \\
               (\epsilon-\sigma\epsilon)^{-1}(1-\sigma\epsilon-f) & \hbox{on $\mathcal{L}(1-\epsilon,1-\sigma\epsilon)\cap E_1$,} \\
               1 & \hbox{otherwise,}
             \end{array}
           \right.
\end{equation*}
where we denote the set $\mathcal{L}(a,b)=\{x\in M|a<f(x)<b\}$.

Denote the set
\begin{equation*}
\Omega=E_1\cap (B_{\rho}(R)\setminus
B_{\rho}(R-1))\cap(\mathcal{L}(\sigma\epsilon,1-\sigma\epsilon)).
\end{equation*}
Recall the growth estimate for $|\nabla f|$ (Corollary 2.3 of
\cite{LW3}):
\begin{equation*}
    \int_{B_{\rho}(R+1)\backslash B_{\rho}(R)}|\nabla f|^2\le
    Ce^{-2R}
\end{equation*}
and the decay estimate for $f$ ((2.10) in \cite{LW3}):
\begin{equation*}\int_{E_1\cap B_{\rho}(R+1)\backslash E_1\cap
B_{\rho}(R)}\rho (1-f)^2\leq Ce^{-2R}.
\end{equation*}
We have
\begin{equation}\label{eqn-c1}
    \left(\int_{\Omega}|\nabla f|^2\right)^\frac{n-2}{n-1}\le
    Ce^{-\frac{2(n-2)}{n-1}R}
\end{equation}
and with notation $S(R)$ as in the statement of theorem,
\begin{equation}\label{eqn-c2}
    \begin{split}
      \int_{\Omega}\rho^{n-1} &\le (S(R))^{2(n-2)}\int_{\Omega}\rho \\
        & \le (\sigma\epsilon)^{-2}(S(R))^{2(n-2)}\int_{\Omega}\rho(1-f)^2\\
        & \le C(S(R))^{2(n-2)}(\sigma\epsilon)^{-2}e^{-2R}.
    \end{split}
\end{equation}
Hence, by $|\nabla r_{\rho}|(x)=\rho(x)$, (\ref{eqn-c1}) and
(\ref{eqn-c2}), we have
\begin{equation}\label{eqn16}\begin{split}
                               \int_{E_1} |\nabla \psi|^2\chi^2|\nabla f|^\frac{2(n-2)}{n-1} & \le  \int_{\Omega}\rho|\nabla
    f|^{\frac {2(n-2)}{n-1}}\\
                                 &\le\left(\int_{\Omega}|\nabla
    f|^2\right)^{\frac
    {n-2}{n-1}}\left(\int_{\Omega}\rho^{n-1}\right)^\frac 1{n-1}\\
    &\le C(\sigma\epsilon)^{-\frac 2{n-1}}(S(R))^{\frac{2(n-2)}{n-1}}e^{-\frac{2(n-2)}{n-1}R-\frac{2}{n-1}R}
    \\
    &\le C(\sigma\epsilon)^{-\frac
2{n-1}}(S(R))^{\frac{2(n-2)}{n-1}}e^{-2R}.
                             \end{split}
\end{equation}
Note the assumption that the Ricci curvature of $M$ is bounded from
below by $-(n-1)\tau(x)$. Then the local gradient estimate of Cheng
and Yau \cite{CY} (cf \cite{LW2}) implies that there exists a
constant $C_n$ depending on $n$ such that
\begin{equation*}
    |\nabla f|(x)\le C_n(\sup_{y\in
B(x,R)}\sqrt{\tau(y)}+R^{-1})|f(x)|, x\in M
\end{equation*}
for  all $R>0$.

Set  $\bar{\rho}(x)=\frac{1}{2}\rho(x)+\frac{1}{2}(n-2)\tau(x), x\in
M$. Then $\sqrt{\tau}\leq \sqrt{\frac{2}{n-2}\bar{\rho}}$ and
\begin{equation}
    |\nabla f|(x)\le
C(\sup_{y\in B(x,R)}\sqrt{\bar{\rho}}(y)+R^{-1})|f(x)|.
\end{equation}

Fix $x\in M$ and consider the function $\eta(R)=\sqrt{2}R-(\sup_{
B(x,R)}\sqrt{\bar{\rho}})^{-1}$. Observe that $\eta(R)$ tends to
$+\infty$ as $R\rightarrow \infty$ and tends to a negative number as
$R\rightarrow 0$. There exists  a $R_0$ depending on $x$ such that
$\sqrt{2}R_0=(\sup_{B(x,R_0)}\sqrt{\bar{\rho}})^{-1}$. hence
\begin{equation}
    |\nabla f|(x)\le C(\sup_{B(x,R_0)}\sqrt{\bar{\rho}})|f(x)|.
\end{equation}

For any $y\in B(x,R_0)$, let $\gamma(s), s\in [0,l]$  be a
minimizing geodesic connecting $x$ and $y$ with respect to the
background metric $ds^2$, where $s$ is the arc-length of $\gamma$ in
$ds^2$. The distance $d_\rho(x,y)$ with respect to $\rho ds^2$
satisfies
\begin{equation}\label{eq15}
    \begin{split}
                   d_{\rho}(x,y)
                     & \leq
\int_{0}^l\sqrt{\rho(\gamma(s))}ds\\
& \le \int_{0}^l\sqrt{2}\sqrt{\bar{\rho}(\gamma(s))}ds\\ &\le
(\sup_{B(x,R_0)}\sqrt{\bar{\rho}})(\sqrt{2}R_0)=1.
                 \end{split}
\end{equation}
This implies  $B(x,R_0)\subset B_\rho(x,1)$. Hence
\begin{equation}
    |\nabla f|(x)\le C(\sup_{B_\rho(x,1)}\sqrt{\bar{\rho}})|f(x)|, x\in M.
\end{equation}
Similarly, we have
\begin{equation}
    |\nabla f|(x)\le C(\sup_{B_\rho(x,1)}\sqrt{\bar{\rho}})|1-f(x)|, x\in M.
\end{equation}
On $E_1$, we have
 \begin{align}\label{e1}
&\int_{E_1} |\nabla \chi|^2\psi^2|\nabla
f|^\frac{2(n-2)}{n-1}\notag\\
&\le C((1-\sigma)\epsilon)^{-2}\int_{\mathcal{L}(1-\epsilon,1-\sigma\epsilon)\cap E_1\cap B_\rho (R)}|\nabla f|^{\frac{2(n-2)}{n-1}+2}\\
&\le
CS^{\frac{2(n-2)}{n-1}}(R+1)((1-\sigma)\epsilon)^{-2}\int_{\mathcal{L}(1-\epsilon,1-\sigma\epsilon)\cap
E_1\cap B_\rho (R)}|\nabla f|^{2}(1-f)^{\frac{2(n-2)}{n-1}}.\notag
\end{align}
Note that Lemma 5.1 of \cite{LW3} asserts that the integral of
$|\nabla f|$ on the level set $l(t)=\{x\in M|f(x)=t\}$, $0\leq t\leq
b$, is invariant. Using this conclusion and the co-area formula and
Lemma 5.1 in \cite{LW3}, we have

\begin{align}\label{esti1}
  &\int_{\mathcal{L}(1-\epsilon,1-\sigma\epsilon)\cap E_1\cap B_\rho (R)}|\nabla
f|^{2}(1-f)^{\frac{2(n-2)}{n-1}}\notag\\
&\qquad\le \int_{1-\epsilon}^{1-\sigma\epsilon}(1-t)^{\frac{2(n-2)}{n-1}}\int_{l(t)\cap E_1\cap B_\rho (R)}|\nabla f|dAdt\notag \\
  &\qquad\le C\int_{l(b)} |\nabla f|
dA\int_{1-\epsilon}^{1-\sigma\epsilon}(1-t)^{\frac{2(n-2)}{n-1}}dt\\
   &\qquad = C\int_{l(b)} |\nabla f|
   dA\cdot
(1-\sigma^{\frac{2(n-2)}{n-1}+1})\epsilon^{\frac{2(n-2)}{n-1}+1}.\notag
\end{align}
Substitute (\ref{esti1}) into (\ref{e1}). Then
\begin{align}\int_{E_1} |\nabla
\chi|^2\psi^2|\nabla f|^\frac{2(n-2)}{n-1}\leq
CS^{\frac{2(n-2)}{n-1}}(R+1)(1-\sigma)^{-2}(1-\sigma^{\frac{2(n-2)}{n-1}+1})\epsilon^{\frac{n-3}{n-1}}.
\end{align}
Setting $\sigma=\frac12$, we have
 $$\int_{E_1}|\nabla \phi|^2|\nabla
f|^{\frac{2(n-2)}{n-1}}\leq
CS^{\frac{2(n-2)}{n-1}}(R+1)(e^{-2R}\epsilon^{-\frac
2{n-1}}+\epsilon^{\frac{n-3}{n-1}}).$$ Let us choose
$\epsilon=e^{-2R}$. Then
\begin{equation}\label{eqn17}\begin{split}
                               \int_{E_1}|\nabla \phi|^2|\nabla f|^{\frac{2(n-2)}{n-1}} \le
CS^{\frac{2(n-2)}{n-1}}(R+1) e^{-\frac{2(n-3)}{n-1}R}.
                             \end{split}
\end{equation}
Using $f$ instead of $1-f$, similar to the above argument, we have
that on $M\backslash E_1$,
\begin{equation}  \int_{M\backslash E_1}|\nabla \phi|^2|\nabla f|^{\frac{2(n-2)}{n-1}} \le CS^{\frac{2(n-2)}{n-1}}(R+1)e^{-\frac{2(n-3)}{n-1})R}.
\end{equation}
Hence
\begin{equation}\label{eqn-m1}  \int_{M}|\nabla \phi|^2|\nabla f|^{\frac{2(n-2)}{n-1}} \le CS^{\frac{2(n-2)}{n-1}}(R+1)e^{-\frac{2(n-3)}{n-1})R}.
\end{equation}
Let $R\to +\infty$. By the assumption on $S(R)$, the left in
(\ref{eqn-m1}) is identically zero. By (\ref{eqn13}), we conclude
that (\ref{eq7}) is actually an equality and hence the improved
Bochner inequality (\ref{boch}) must be an equality. Note that Lemma
4.1 of \cite{LW3} asserts that if equality in inequality
(\ref{boch}) holds, the metric of $M$ must be a warped product as
described in the theorem. We obtain the conclusion of theorem in the
case of $n\geq 4$.

 In the case of $n=3$, we may choose $\psi$ as above and $\chi$ to be
\begin{equation*}
    \chi(x)=\left\{
             \begin{array}{ll}
               0 & \hbox{ on $\mathcal{L}(0,\sigma\epsilon)\cup \mathcal{L}(1-\sigma\epsilon,1)$,} \\
               (-\log\sigma)^{-1}(\log f-\log(\sigma\epsilon))  & \hbox{on $\mathcal{L}(\sigma\epsilon,\epsilon)\cap(M\setminus E_1)$,} \\
               (-\log\sigma)^{-1}(\log(1-f)-\log(1-\sigma\epsilon)) & \hbox{on $\mathcal{L}(1-\epsilon,1-\sigma\epsilon)\cap E_1$,} \\
               1 & \hbox{otherwise.}
             \end{array}
           \right.
\end{equation*}

By an argument similar to the above one for $n\geq 4$ (combining
with the corresponding estimates for $n=3$ in Theorem 5.2 in
\cite{LW3}), we have the estimate

\begin{equation}\label{esti3}\int_{M}|\nabla \phi|^2|\nabla f|^{\frac{2(n-2)}{n-1}} \le CS(R+1)(\sigma^{-1}\epsilon^{-1}e^{-2R}+(-\log \sigma)^{-1}).
\end{equation}

Choose $\sigma=\epsilon=e^{-Rq(R)}$ with
$q(R)=\sqrt{\frac{S(R+1)}{R}}$. Then using the argument in
\cite{LW3}, we have the right side of (\ref{esti3}) tends to zero as
$R\rightarrow +\infty$. We conclude that (\ref{eq7}) is actually an
equality and hence the theorem holds for $n=3$.

\end{proof}

\section{Application  to minimal hypersurfaces}

Let $M^n$ be a complete minimal minimal hypersurface in
$\mathbb{R}^{n+1}$ for $ n\geq 3$. We first give some examples of
the metric $\rho ds^2$ such that $M$ satisfies property
$(\mathcal{P}_{\rho})$.

{\bf Example 3.1.} Let $\bar{d}(x,y), x,y \in \mathbb{R}^{n+1}$
denote the distance between $x$ and $y$ in $\mathbb{R}^{n+1}$.
Denote by $\bar{r}(x), x\in M$  the extrinsic distance function
$\bar{d}(x,o)$ from a fixed point $o\in \mathbb{R}^3$ ($o$ may be in
$M$ or not). It is known that
$$\Delta\bar{r}\geq (n-1)\bar{r}^{-1},$$where $\Delta$ is the
Laplacian on $M$.

For any $\phi\in C^{\infty}_o(M)$,
\begin{align*}
&(n-1)\int_M\bar{r}^{-2}\phi^2\leq\int_M\bar{r}^{-1}\phi^2\Delta\bar{r}\\
&=-2\int_M\bar{r}^{-1}\phi\langle\nabla\phi,\nabla\bar{r}\rangle+\int_M\bar{r}^{-2}\phi^2|\nabla\bar{r}|^2\\
&\leq 2\int_M\bar{r}^{-1}\phi|\nabla\phi|+\int_M\bar{r}^{-2}\phi^2.
\end{align*}
\begin{align*}(n-2)\int_M\bar{r}^{-2}\phi^2&\leq2\int_M\bar{r}^{-1}\phi|\nabla\phi|\\
&\leq2(\int_M\bar{r}^{-2}\phi^2)^{\frac12}(\int_M|\nabla\phi|^2)^{\frac12}
\end{align*}
Hence
\begin{equation}\label{eqn-p}\int_M|\nabla \phi|^2\ge \frac{(n-2)^2}{4}\int_M \bar{r}^{-2}\phi^2 \quad \text{for all}\quad  \phi\in
C_o^{+\infty}(M).
\end{equation}
Let $\rho(x)=\frac{(n-2)^2}{4}\bar{r}^{-2}(x), x\in M$. Inequality
(\ref{eqn-p}) asserts the Poincar\'e inequality holds with weighted
function $\rho$.

 Further the metric
$\rho ds^2$ is complete. Indeed, take a fixed point $p\in M$ with
$p\neq o$. Let $r(x), x\in M$ denote the intrinsic distance from
$p$. Note that $\bar{r}(x)\leq\bar{d}(o,p)+\bar{d}(x,p)\leq
r_0+r(x)$, where $r_0=\bar{d}(o,p)>0$. Then
$\bar{r}^{-2}(x)>(r_0+r(x))^{-2}.$ It is known that the metric
$(r_0+r(x))^{-2}ds^2$ is complete. Hence $\rho ds^2$ is complete.

Thus we obtain that $M$ has  property $(\mathcal{P}_\rho)$ for
$\rho$.

{\bf Example 3.2.} Using   smoothing technique, we may modify
$\rho=\frac{(n-2)^2}{4}\bar{r}^{-2}$ in Example 2.1 to get a bounded
smooth positive function $\rho_1(x)=\rho_1(\bar{r}(x))$, $x\in M$
such that M has property $(\mathcal{P}_\rho)$ for  $\rho_1$.

Indeed, let positive number $0<b\leq r_0$ fixed, we can choose
number $a, 0<a<b$ such that function
$\zeta(\bar{r})=\frac{(n-2)^2}{4}(\bar{r}^{-2}-e^{-\frac{1}{(\bar{r}-b)^2}})$
is strictly decreasing in $(a,b)$ as $\bar{r}$ tends increasingly to
$b$ and construct the smooth $\rho_1$
\begin{equation*}
\rho_1(\bar{r}(x))=\left\{
                               \begin{array}{ll}
                                  h(\bar{r}) & \hbox{for $\bar{r}(x)\leq a$,} \\
                                 \zeta(\bar{r})  & \hbox{for $a<\bar{r}(x)<b$,} \\
                                 \rho(\bar{r}) & \hbox{for $\bar{r}(x)\geq b$,}
                               \end{array}
                             \right.
\end{equation*}
where $h(\bar{r})$ is chosen to be bounded and to satisfy $
\rho(\bar{r})\geq h(\bar{r})\geq \rho(b)$ for $\bar{r}\leq a$.

Observe that $\rho_1\leq \rho$. Hence the Poincar\'e inequality
holds for $\rho_1$. Moreover $\rho_1ds^2$ is complete since
$\rho_1(x)\geq \frac{(n-2)^2}{4}(r_0+r(x))^{-2}$. In fact, for
$\bar{r}(x)\geq b$, $\rho_1=\rho$.  Note that for $\bar{r}(x)< b$,
$\rho_1(\bar{r}(x)\geq \rho(b)$ and $0<b\leq r_0$. Hence
$\rho_1(\bar{r}(x))\geq\frac{(n-2)^2}{4}(r_0+r(x))^{-2}$ for
$\bar{r}(x)< b$.

{\bf Example 3.3} Under the above notations, choose
$\rho_2(x)=\frac{(n-2)^2}{4}(r_0+r(x))^{-2}, x\in M$. Since
$\rho_2\leq \rho$, Poincar\'e inequality holds with weighted
function $\rho_2$. By the completeness of the metric $\rho_2ds^2$,
we know $M$ has property $(\mathcal{P}_\rho)$ for $\rho_2$.

\begin{thm}\label{thm-01-1}(Theorem \ref{thm01})
Let $M$ be an $\frac{n-2}{n}$-stable complete minimal hypersurface
in $\mathbb{R}^{n+1}, n\geq 3$ and the norm of its
 second fundamental form satisfies
\begin{align}
    &\lim_{R\to +\infty}\sup_{B(R)}|A|
/R^{\frac{n-3}{2}}=0 \quad \text{for}\quad n\geq 4; \notag\\
&\lim_{R\to +\infty}\sup_{B(R)}|A|/{\ln R}=0, \quad \text{for}\quad
n=3,
\end{align} then $M$ either has one end or must be a
catenoid.
\end{thm}

\begin{proof}By the Gauss equation,\begin{equation}\label{eq5}
    \begin{split}
      |A|^2 & \ge h_{11}^2+\sum_{i=2}^nh_{ii}^2+2\sum_{i=1}^nh_{1i}^2\\
        &\ge h_{11}^2+\frac{(\sum_{i=2}^nh_{ii})^2}{n-1}+2\sum_{i=1}^nh_{1i}^2\\
        &\ge \frac{n}{n-1}\left(h_{11}^2+\sum_{i=2}^nh_{1i}^2\right)\\
        &\ge -\frac{n}{n-1}Ric_M(e_1,e_1).
    \end{split}
    \end{equation}
Let us choose $\tau=\frac{|A|}{n}$ and $\rho=\rho_1$ (or $\rho_2$)
in Theorem \ref{thm-r01}. By the boundedness of $\rho_1$ (or
$\rho_2$), the growth assumption (\ref{eq}) on $\rho$ is satisfied.
Now we will assert that the growth assumption (\ref{eq}) on $\tau$
is also satisfied.

It  can be directly verified that a minimizing geodesic starting
from the fixed point $p$ with respect to $ds^2$ is also a minimizing
geodesic starting from $p$ with respect to $\rho_2ds^2$. Then by
direct calculation,  we have $B_{\rho_2}(\bar{R})=B(R)$, where
$\bar{R}=\frac{n-2}{2}\ln(1+\frac{R}{r_0})$. Then for $n\geq 4$
\begin{align}
\lim_{\bar{R}\to +\infty}\sup_{B_{\rho_2}(\bar{R})}|A|
e^{-\frac{(n-3)}{n-2}\bar{R}}=C\lim_{R\to
+\infty}\sup_{B(R)}|A|R^{-\frac{n-3}{2}}=0.
\end{align}
For $n=3$,
\begin{equation}\lim_{\bar{R}\to +\infty}\sup_{B_{\rho_2}(\bar{R})}|A|
\bar{R}^{-1}=C\lim_{R\to +\infty}\sup_{B(R)}|A|(\ln R)^{-1}=0.
\end{equation}

If $\rho=\rho_1$, by $\rho_1\geq\rho_2$,
$B_{\rho_1}(\bar{R})\subseteq B_{\rho_2}(\bar{R})$ and hence the
growth assumption on $\tau$  also holds for $\rho_1$.

Therefore the conclusion of Theorem \ref{thm-r01} is valid. Let us
assume that $M$ has at least two ends. Since every end of a complete
noncompact minimal hypersurface in $\mathbb{R}^{n+1}$ is
nonparabolic (\cite{CSZ}, see its proof also in \cite{CCZ}), by
Theorem \ref{thm-r01}, we know that $M$ has exactly two nonparabolic
ends and $M=\mathbb{R}\times N$ with the warped product metric
$$ds_M^2=dt^2+\eta^2(t)ds^2_N,$$
for some compact manifold $N$ and some positive function $\eta(t)$.
Moreover, $|A|$ is a function of $t$ alone satisfying
$$(n-2)\eta''\eta^{-1}=\frac{|A|}{n}.$$ Hence $M$
has a rotationally symmetric metric. By a result of do Carmo and
Dajczer (\cite{dcd}, Corollary 4.4), it implies that every part of M
is a part of a catenoid. Hence $M$ is contained in a catenoid
$\mathcal{C}$ by minimality of the immersion. Since $M$ is complete
and  the catenoid $\mathcal{C}$ is simply connected because $n\ge3$,
$M$ must be an embedded hypersurface, see \cite[p.330]{Sp}. Hence
$M$ is the catenoid.
\end{proof}
Theorem \ref{thm-01-1} implies directly that

\begin{cor}(Corollary \ref{bound}) Let $M$ be an $\frac{n-2}{n}$-stable complete minimal
hypersurface in $\mathbb{R}^{n+1}, n\geq 3$ with at least two ends.
If its
 second fundamental form is bounded, then $M$ must be a
catenoid.
\end{cor}

\bigskip
\bigskip
\noindent  Xu Cheng\\Insitituto de Matem\'atica\\Universidade
Federal Fluminense - UFF\\Centro, Niter\'{o}i, RJ 24020-140 Brazil
\\e-mail:xcheng@impa.br

\bigskip
\bigskip
\noindent Detang Zhou\\Insitituto de Matem\'atica\\Universidade
Federal Fluminense - UFF\\Centro, Niter\'{o}i, RJ 24020-140 Brazil
\\e-mail: zhou@impa.br
\end{document}